\documentclass[11pt]{article}

\usepackage[T1]{fontenc}
\usepackage[utf8]{inputenc}
\usepackage{lmodern}
\usepackage[a4paper,margin=27mm]{geometry}
\usepackage{amsmath,amssymb,amsthm}
\usepackage{booktabs,longtable,array}
\usepackage{microtype}
\usepackage[hidelinks]{hyperref}
\usepackage{enumitem}
\usepackage{tikz}

\setlist{itemsep=3pt,topsep=5pt}
\newtheorem{theorem}{Theorem}[section]
\newtheorem{lemma}[theorem]{Lemma}
\newtheorem{proposition}[theorem]{Proposition}
\newtheorem{corollary}[theorem]{Corollary}

\theoremstyle{definition}

\theoremstyle{remark}

\newcommand{\RCB}{\operatorname{RCB}}
\newcommand{\brank}{\operatorname{rank}_{\mathbb B}}
\newcommand{\supp}{\operatorname{supp}}
\newcommand{\A}{\mathcal A}
\newcommand{\B}{\mathcal B}
\newcommand{\C}{\mathcal C}
\newcommand{\Om}{\Omega}
\newcommand{\code}[1]{\texttt{\detokenize{#1}}}

\title{Boolean Rank Nine for Polygons with 35 to 40 Vertices}

\author{
Dake Zhu\\[4pt]
{\small Qiuzhen College, Tsinghua University}\\
{\small Beijing, China}
}

\date{1 November 2026}

\begin{document}

\maketitle

\begin{abstract}
We give a computer-assisted proof that the rectangle covering number,
or Boolean rank, of every polygon slack matrix of order 35 through 40
is nine. The upper bound is due to Goucha, Gouveia, and Silva; the new
contribution is the matching lower bound. A hypothetical covering by
eight rectangles yields a saturated inclusion cycle between two
antichains on eight coordinates. An exact maximal-chain identity bounds
its length and assigns a nonnegative penalty to each inclusion edge.
The equality case is excluded combinatorially, leaving at most 39 lower
vertices. For orders 35 through 39, symmetry reductions and a shadow
bound leave fifteen finite exclusions, each supported by a checked
reverse-unit-propagation refutation. The main text presents the
mathematical reductions and their dependencies; the accompanying
supplement gives the complete cycle parameters, encodings, and
certificate record. The previously known equality for orders 41 through
55 is cited directly.
\end{abstract}

\medskip


\begingroup
\small
\tableofcontents
\endgroup

\medskip
\newpage
\section{Problem, prior results, and proof strategy}
\label{sec:introduction}

\subsection{Polygon supports and rectangle coverings}

The support of a polygon slack matrix has a simple cyclic form, even
though its positive entries depend on the geometry. Let $P$ be a convex
polygon with exactly $n$ vertices and $n$ irredundant facet inequalities.
After cyclically indexing rows and columns, its slack matrix $S_P$ has
the support of the zero--one matrix $M_n$ defined by
\begin{equation}
(M_n)_{ij}=0
\quad\Longleftrightarrow\quad
i=j\ \text{or}\ i=j+1\pmod n.
\label{eq:support}
\end{equation}
All other entries are one.
Explicitly, the support matrix has the cyclic form
\begin{equation*}
M_n=
\begin{pmatrix}
0 & 1 & 1 & \cdots & 1 & 0\\
0 & 0 & 1 & \cdots & 1 & 1\\
1 & 0 & 0 & \cdots & 1 & 1\\
\vdots & \ddots & \ddots & \ddots & \vdots & \vdots\\
1 & \cdots & 1 & 0 & 0 & 1\\
1 & \cdots & 1 & 1 & 0 & 0
\end{pmatrix}.
\label{eq:Mn-matrix}
\end{equation*}
A rectangle of a nonnegative matrix $M$ is a set $I\times J$ contained
in $\supp(M)$. The rectangle covering number $\RCB(M)$ is the smallest
number of such rectangles whose union is $\supp(M)$. It equals the
Boolean rank $\brank(M)$: the smallest inner dimension of a Boolean
factorization of the support, with AND for multiplication and OR for
addition.

For a nonnegative matrix, we use $\brank(M)$ to mean the Boolean rank
of its zero--one support matrix. Thus
\[
\RCB(S_P)=\brank(S_P)=\RCB(M_n)=\brank(M_n).
\]

\subsection{Published bounds and the remaining gap}

Goucha, Gouveia, and Silva~\cite{GGS} studied polygon Boolean ranks
through homogeneous factorizations and cycles in Johnson graphs.
Their Section 3.2, Table 1 gives a size-nine factorization for every
$35\le n\le55$. In the interval $35\le n\le40$, the asterisk in that
table marks an upper bound whose optimality was not established.
For $41\le n\le55$, the equality with nine is already established.
Engel, Kalinowski, and Schymura~\cite[Section 5.7, Table 3]{EKS}
also list the bounds $8\le\RCB(M_n)\le9$ for $35\le n\le40$.

\begin{proposition}[Published upper bound]
\label{prop:known-upper}
For every $35\le n\le55$, one has $\RCB(M_n)\le9$.
\end{proposition}

\begin{proof}
Goucha, Gouveia, and Silva give homogeneous Boolean factorizations of
size nine throughout this interval~\cite[Section 3.2, Table 1]{GGS}.
Every homogeneous Boolean factorization is also an unrestricted Boolean
factorization, so it gives a covering by at most nine rectangles.
\end{proof}

The published size-nine witnesses and their entrywise checks are
documented in the accompanying \emph{Supplementary Material}, Section S5.

\subsection{Main result}

\begin{theorem}
\label{thm:uniform}
For every convex polygon $P$ with $35\le n\le40$ vertices,
\[
\RCB(S_P)=\brank(S_P)=9.
\]
\end{theorem}

\begin{corollary}
\label{cor:range}
For every convex polygon $P$ with $35\le n\le55$ vertices,
\[
\RCB(S_P)=\brank(S_P)=9.
\]
\end{corollary}

\begin{proof}
Use Theorem~\ref{thm:uniform} for $35\le n\le40$.
For $41\le n\le55$, use the exact values already established in
\cite[Section 3.2, Table 1]{GGS}.
\end{proof}

\subsection{Proof strategy and dependencies}
\label{sec:strategy}

Only the lower bound for $35\le n\le40$ is new here. Since the upper
bound is already available, the task is to exclude a covering by at most
eight rectangles. The proof has four steps:
\begin{enumerate}[label=\textup{(\arabic*)}]
\item
Section~\ref{sec:normal} turns such a covering into a saturated
inclusion cycle on eight coordinates.

\item
Section~\ref{sec:chains} proves rank restrictions and an exact
penalty budget by maximal-chain counting.

\item
Section~\ref{sec:analytic} excludes equality in that budget,
giving $N\le39$ without a computational refutation.

\item
Section~\ref{sec:reductions} excludes the remaining cycles through
the H, G, T, M, and F stages. Section~\ref{sec:main-proof} then assembles
the argument.
\end{enumerate}

Here $n$ denotes the polygon order and $N$ the number of lower vertices
in an inclusion cycle; the normal form gives $N=n$. The saved selection
formulas allow every $35\le N\le40$ simultaneously. Thus the finite
exclusions use no assumption of homogeneity or monotonicity in $n$.

The local stages are applied in a fixed order. H excludes homogeneous
cycles. G restricts the possible rank pairs. T removes the remaining
lower two-sets and, by duality, upper six-sets. M proves that mixed edges
form a matching. Only then does F use this matching condition.
The complete encodings and certificates are supplied in the supplement.

\section{From rectangle coverings to inclusion cycles}
\label{sec:normal}

Write $\Om_k=\{0,1,\ldots,k-1\}$ and use cyclic indices in
$\mathbb Z/n\mathbb Z$. The \emph{rank} of a set means its cardinality.
The sets below record which rectangles contain a given row or column;
their ground-set coordinates label rectangles, not polygon vertices.
If a covering has fewer than eight rectangles, pad it with empty
rectangles. It is therefore enough to exclude coverings indexed by
$\Om_8$.

\subsection{The maximal-rectangle normal form}

\begin{lemma}[Maximal-rectangle normal form]
\label{lem:normal}
Suppose $M_n$ has a covering by $k$ rectangles:$I_{0}\times J_{0},I_{1}\times J_{1},\dots ,I_{k-1}\times J_{k-1}$. There is a covering by
$k$ rectangles $I_{0}^{*}\times J_{0}^{*},I_{1}^{*}\times J_{1}^{*},\dots ,I_{k-1}^{*}\times J_{k-1}^{*}$, with row membership sets $A_i=\{t\in \Om_k : i\in I_{t}^{*}\}$ and column
membership sets $D_j=\{t\in \Om_k : j\in J_{t}^{*}\}$. Let
$B_j=\Om_k\setminus D_j$, then we have
\begin{align}
B_j&=A_j\cup A_{j+1},
\label{eq:union}\\
A_i&=B_{i-1}\cap B_i,
\label{eq:intersection}\\
A_i\subseteq B_j
&\quad\Longleftrightarrow\quad
i\in\{j,j+1\}.
\label{eq:crown}
\end{align}
\end{lemma}

\begin{proof}
For a rectangle $I\times J$, let $J^*$ consist of every column
compatible with all rows in $I$, and let $I^*$ consist of every row
compatible with all columns in $J^*$. Then $I\subseteq I^*$ and
$J\subseteq J^*$. A column outside $J^*$ has a zero in some row of
$I$, hence also in a row of $I^*$. Thus $I^*\times J^*$ is maximal
on both sides. Extend every rectangle this way. No previously covered
positive position is lost.

For a maximal rectangle indexed by $t$, column $j$ belongs to the
rectangle exactly when neither row $j$ nor row $j+1$ belongs to it.
Consequently $t\in D_j$ exactly when
$t\notin A_j\cup A_{j+1}$, proving~\eqref{eq:union}. Similarly, row $i$
belongs to the rectangle exactly when neither column $i-1$ nor column
$i$ belongs to it. This gives~\eqref{eq:intersection}.

Finally, position $(i,j)$ is uncovered exactly when
$A_i\cap D_j=\varnothing$, or equivalently $A_i\subseteq B_j$.
The prescribed zero positions are precisely those in~\eqref{eq:crown}.
\end{proof}

\subsection{Antichains and the induced cycle}

\begin{lemma}
\label{lem:antichains}
For $n\ge4$, sets satisfying~\eqref{eq:crown} form two antichains
$\A=\{A_i\}$ and $\B=\{B_j\}$, each of cardinality $n$. No member of
$\A$ equals a member of $\B$. The bipartite graph defined by strict
inclusion $A\subsetneq B$ is an induced alternating cycle of length
$2n$.
\end{lemma}

\begin{proof}
If $A_i\subseteq A_h$, then $A_i$ lies in both $B_{h-1}$ and $B_h$.
By~\eqref{eq:crown}, $i$ belongs to both $\{h-1,h\}$ and
$\{h,h+1\}$, so $i=h$. This proves the first antichain assertion and
also distinctness of the $A_i$.

If $B_j\subseteq B_h$, its two lower neighbors are also lower neighbors
of $B_h$. Equation~\eqref{eq:crown} then implies
$\{j,j+1\}=\{h,h+1\}$ and hence $j=h$. If $A_i=B_j$, the other
lower neighbor of $B_j$ would be contained in $A_i$, contrary to the
first assertion. The last statement now follows directly
from~\eqref{eq:crown}.
\end{proof}

Conversely, a cyclic sequence satisfying~\eqref{eq:crown}, with
$B_j=A_j\cup A_{j+1}$, yields a covering. Its rectangle indexed by $t$ is
\begin{equation}
\{i:t\in A_i\}
\times
\{j:t\notin A_j\cup A_{j+1}\}.
\label{eq:rectangles}
\end{equation}
We call a pair of antichains satisfying
\eqref{eq:union}--\eqref{eq:crown} a \emph{saturated inclusion cycle}.
This is a condition on the selected antichains, not a claim that the
cycle is maximal among cycles in the entire Boolean lattice.

\begin{figure}[htbp]
\centering
\begin{tikzpicture}[
    font=\small,
    Avertex/.style={
        circle,
        draw,
        fill=white,
        inner sep=0pt,
        minimum size=5pt
    },
    Bvertex/.style={
        rectangle,
        draw,
        fill=black!25,
        inner sep=0pt,
        minimum size=5pt
    },
    every label/.style={font=\scriptsize},
    line width=0.6pt
]

\foreach \i in {0,...,5} {
    \node[Avertex, label=below:{$A_{\i}$}]
        (LA\i) at ({0.88*\i},0) {};
    \node[Bvertex, label=above:{$B_{\i}$}]
        (LB\i) at ({0.88*\i},1.8) {};
}

\foreach \i in {0,...,5} {
    \draw (LA\i) -- (LB\i);
}

\foreach \i/\j in {0/1,1/2,2/3,3/4,4/5} {
    \draw (LB\i) -- (LA\j);
}

\draw (LB5)
    -- (4.9,1.8)
    -- (4.9,2.55)
    -- (-0.45,2.55)
    -- (-0.45,0)
    -- (LA0);

\node at (2.2,-1.15)
    {(a) Bipartite inclusion graph};

\draw[->, >=stealth]
    (5.35,1.05) -- (6.35,1.05)
    node[midway,above,font=\scriptsize] {redraw};

\begin{scope}[shift={(8.65,1.05)}]

    \foreach \i in {0,...,5} {
        \pgfmathtruncatemacro{\angleA}{90-60*\i}
        \pgfmathtruncatemacro{\angleB}{60-60*\i}

        \node[Avertex,
              label={\angleA:{$A_{\i}$}}]
            (RA\i) at (\angleA:1.45) {};

        \node[Bvertex,
              label={\angleB:{$B_{\i}$}}]
            (RB\i) at (\angleB:1.45) {};
    }

    \foreach \i in {0,...,5} {
        \pgfmathtruncatemacro{\nextindex}{mod(\i+1,6)}
        \draw (RA\i) -- (RB\i) -- (RA\nextindex);
    }

    \node at (0,0) {$C_{12}$};

\end{scope}

\node at (8.65,-1.15)
    {(b) Alternating cycle};

\end{tikzpicture}
\caption{
$n=6$.
}
\label{fig:inclusion-cycle}
\end{figure}

\subsection{Complement duality and homogeneity}

\begin{lemma}[Complement duality]
\label{lem:dual}
Taking $\tilde{A_i}=B_{-i}^{c}$, $\tilde{B_i}= A_{-i}^{c}$ preserves a
saturated inclusion cycle. An inclusion edge $A\subsetneq B$ has rank type $(|A|,|B|)$;
under this transformation, rank type $(a,b)$ becomes
$(k-b,k-a)$.

\end{lemma}

\begin{proof}
Taking complements and exchanging the two families preserves
the antichain property. Moreover,
\begin{align*}
\widetilde A_i\cup\widetilde A_{i+1}
&=(B_{-i}\cap B_{-i-1})^c
 =A_{-i}^c
 =\widetilde B_i,\\
\widetilde B_{i-1}\cap\widetilde B_i
&=(A_{1-i}\cup A_{-i})^c
 =B_{-i}^c
 =\widetilde A_i.
\end{align*}
The inclusion condition is also preserved:
\[
\widetilde A_i\subseteq\widetilde B_j
\quad\Longleftrightarrow\quad
A_{-j}\subseteq B_{-i}
\quad\Longleftrightarrow\quad
i\in\{j,j+1\}.
\]
Finally, an edge $A\subsetneq B$ becomes
$B^c\subsetneq A^c$, whose rank type is $(k-b,k-a)$.
\end{proof}

A cycle is \emph{homogeneous of type $(r,r+1)$} if every lower set
has size $r$ and every upper set has size $r+1$.

\section{Maximal-chain counting and the penalty budget}
\label{sec:chains}

\subsection{The exact chain identity}

A \emph{chain} in $2^{\Omega_k}$ is a family of subsets ordered
by inclusion. It is \emph{maximal} if no additional subset can be
inserted while preserving this ordering. Every maximal chain has
the form
\[
\mathcal C=\{C_0,C_1,\ldots,C_k\},
\qquad
\varnothing=C_0\subsetneq C_1\subsetneq\cdots
\subsetneq C_k=\Omega_k,
\qquad |C_r|=r.
\]
Thus a maximal chain contains exactly one set of each cardinality.
It is uniquely determined by a permutation
$(\pi_1,\ldots,\pi_k)$ of $\Omega_k$, through
\[
C_0=\varnothing,
\qquad
C_r=\{\pi_1,\ldots,\pi_r\}
\quad (1\le r\le k).
\]
Choosing the permutation uniformly therefore gives a uniformly
random maximal chain.

We say that $\mathcal C$ \emph{meets} a family $\mathcal A$
if some set on the chain belongs to $\mathcal A$; equivalently,
\[
\mathcal C\cap\mathcal A\ne\varnothing.
\]
If $\mathcal A$ is an antichain, the chain can contain at most
one member of $\mathcal A$, since any two distinct sets on the
chain are comparable by inclusion.

\begin{lemma}[Chain identity]
\label{lem:chain}
Let $\A,\B$ be disjoint antichains in $2^{\Om_k}$, with every
cross-comparability oriented from $\A$ to $\B$. Let $E$ be the set of
pairs $(A,B)$ with $A\subsetneq B$. A uniformly random maximal chain
$\C$ satisfies
\begin{align}
\Pr\bigl(\C\cap(\A\cup\B)\ne\varnothing\bigr)
&=
\sum_{A\in\A}\binom{k}{|A|}^{-1}+
\sum_{B\in\B}\binom{k}{|B|}^{-1}-
\sum_{(A,B)\in E}
\frac{1}{\binom{k}{|B|}\binom{|B|}{|A|}}
\le1.
\label{eq:chain-identity}
\end{align}
\end{lemma}

\begin{proof}
There are $k!$ maximal chains. A set of size $a$ lies on
$a!(k-a)!$ of them. A nested pair of sizes $a<b$ lies on
$a!(b-a)!(k-b)!$ of them. A chain meets each antichain at most once.
If it meets both, the two selected sets form a pair in $E$. Thus
inclusion--exclusion stops after the pair terms and gives the displayed
identity.
\end{proof}

In the cycles and paths considered below, every upper vertex
has two lower neighbors. This excludes reverse
cross-comparabilities. Indeed, if $B\subseteq A$ for some
$B\in\mathcal B$ and $A\in\mathcal A$, a lower neighbor
$A'$ of $B$ would satisfy
\[
A'\subsetneq B\subseteq A,
\]
contradicting the antichain property of $\mathcal A$.

Define, for $1\le a<b\le k-1$,
\begin{equation}
w_k(a,b)
=
\frac{1}{2\binom{k}{a}}
+
\frac{1}{2\binom{k}{b}}
-
\frac{1}{\binom{k}{b}\binom{b}{a}}.
\label{eq:weight}
\end{equation}
In particular, the eight-coordinate formula is
\begin{equation}
\boxed{\displaystyle
w_8(a,b)
=
\frac{1}{2\binom{8}{a}}
+
\frac{1}{2\binom{8}{b}}
-
\frac{1}{\binom{8}{b}\binom{b}{a}}.
}
\label{eq:weight-eight}
\end{equation}
For an inclusion cycle, every selected vertex has degree two. Therefore
Lemma~\ref{lem:chain} becomes
\begin{equation}
\sum_{(A,B)\in E}w_k(|A|,|B|)\le1.
\label{eq:cycle-weight}
\end{equation}
For an induced alternating path with $m$ lower vertices and $m-1$ upper
vertices, whose endpoints are $A_L,A_R$, the corresponding inequality is
\begin{equation}
\sum_{(A,B)\in E}w_k(|A|,|B|)
+
\frac{1}{2\binom{k}{|A_L|}}
+
\frac{1}{2\binom{k}{|A_R|}}
\le1.
\label{eq:path-weight}
\end{equation}

\subsection{Analytic bounds for the chain weights}

\begin{lemma}[Analytic weight minima]
\label{lem:weight-minima}
The chain weights satisfy
\begin{align}
\min_{1\le a<b\le6}w_7(a,b)
&=\frac3{140},
\label{eq:min-seven}\\
\min_{\substack{2\le a\le5,\ 3\le b\le6\\a<b}}w_8(a,b)
&=\frac1{80}.
\label{eq:min-eight}
\end{align}
Equality in~\eqref{eq:min-seven} holds only at $(3,4)$;
equality in~\eqref{eq:min-eight} holds only at $(3,4)$ and $(4,5)$.
\end{lemma}

\begin{proof}

For $k=7$, all fifteen admissible pairs are listed below:
\[
\begin{array}{c|rrrrr}
420w_7(a,b) & b=2 & b=3 & b=4 & b=5 & b=6\\
\hline
a=1 & 30 & 32 & 33 & 36 & 50\\
a=2 & \text{--} & 12 & 14 & 18 & 36\\
a=3 & \text{--} & \text{--} & 9 & 14 & 33\\
a=4 & \text{--} & \text{--} & \text{--} & 12 & 32\\
a=5 & \text{--} & \text{--} & \text{--} & \text{--} & 30
\end{array}
\]
The smallest entry is $9$, attained only at $(3,4)$.
Thus the minimum weight is $9/420=3/140$.

For $k=8$, all ten admissible pairs in the stated range are:
\[
\begin{array}{c|rrrr}
1680w_8(a,b) & b=3 & b=4 & b=5 & b=6\\
\hline
a=2 & 35 & 38 & 42 & 56\\
a=3 & \text{--} & 21 & 27 & 42\\
a=4 & \text{--} & \text{--} & 21 & 38\\
a=5 & \text{--} & \text{--} & \text{--} & 35
\end{array}
\]
The smallest entry is $21$, attained only at $(3,4)$ and $(4,5)$.
Thus the minimum weight is $21/1680=1/80$.
\end{proof}

\subsection{The seven-coordinate path bound and rank restrictions}

\begin{lemma}
\label{lem:seven-path}
An induced alternating inclusion path on a seven-element ground set,
with two antichains and both endpoints in the lower antichain, has at
most 23 lower vertices.
\end{lemma}

\begin{proof}
It suffices to rule out $m\ge24$. In that case neither antichain contains
the empty set or the full ground set: these would violate either the
antichain property or the degree bound of two. There are no reverse
cross-comparabilities, by the same antichain argument as above, since
each upper vertex has two lower neighbors. Thus the relevant ranks
satisfy $1\le a<b\le6$.

Lemma~\ref{lem:weight-minima} gives $w_7(a,b)\ge3/140$, while binomial
unimodality gives $\binom{7}{a}\le35$. There are $2(m-1)$ edges, so
\eqref{eq:path-weight} implies
\[
1
\ge
2(m-1)\frac3{140}+\frac1{35}
=
\frac{3m-1}{70}.
\]
Hence $m\le23$.
\end{proof}

\begin{proposition}
\label{prop:ranks}
Every saturated inclusion cycle with $N\ge25$ lower vertices in
$2^{\Om_8}$ satisfies
\begin{equation}
2\le|A|\le5
\quad(A\in\A),
\qquad
3\le|B|\le6
\quad(B\in\B).
\label{eq:rank-range}
\end{equation}
\end{proposition}

\begin{proof}
We first show that $|A|\ge2$ for every $A\in\mathcal A$.
The empty set cannot belong to $\mathcal A$, since $\mathcal A$
is an antichain with more than one member.

Suppose that $A_i=\{t\}$.
Every other lower set omits $t$, by the antichain property.
Delete $A_i$ and its two upper neighbors.
Each remaining upper set is the union of two remaining lower
sets, so it also omits $t$.
The resulting induced alternating path therefore lies on the
seven-element ground set $\Omega_8\setminus\{t\}$ and has
$N-1\ge24$ lower vertices, contradicting
Lemma~\ref{lem:seven-path}.
Thus $|A|\ge2$ for every $A\in\mathcal A$.

Apply this lower bound to the complement-dual cycle from
Lemma~\ref{lem:dual}.
Its lower sets are the complements of the original upper sets.
Consequently,
\[
8-|B|\ge2,
\qquad\text{so}\qquad
|B|\le6
\quad(B\in\mathcal B).
\]

Every lower set $A$ has an upper neighbor $B$ with
$A\subsetneq B$. Hence
\[
|A|\le |B|-1\le5.
\]
Finally, every upper set $B$ has a lower neighbor $A$ with
$A\subsetneq B$. Using the lower bound already proved gives
\[
|B|\ge |A|+1\ge3.
\]
Together, these bounds yield
\[
2\le |A|\le5
\quad(A\in\mathcal A),
\qquad
3\le |B|\le6
\quad(B\in\mathcal B).
\]
\end{proof}

\subsection{The common penalty budget}

For the ranks in~\eqref{eq:rank-range}, define the integer penalty
\begin{equation}
p(a,b)
=
1680\left(w_8(a,b)-\frac1{80}\right).
\label{eq:penalty}
\end{equation}
Lemma~\ref{lem:weight-minima} proves $p(a,b)\ge0$ analytically.
For example,
\[
w_8(3,4)
=
\frac1{112}+\frac1{140}-\frac1{280}
=
\frac1{80}.
\]
Substitution in~\eqref{eq:weight-eight}, using
$w_8(a,b)=w_8(8-b,8-a)$, gives the complete exact penalty table
\begin{equation}
\begin{array}{c|rrrrrrrrrr}
(a,b)
&(2,3)&(2,4)&(2,5)&(2,6)&(3,4)
&(3,5)&(3,6)&(4,5)&(4,6)&(5,6)\\
\hline
p(a,b)
&14&17&21&35&0&6&21&0&17&14
\end{array}
\label{eq:penalty-table}
\end{equation}
These are rational identities, rather than numerical approximations.

There are $2N$ inclusion edges. Subtracting their baseline contribution
$2N/80$ from~\eqref{eq:cycle-weight} and multiplying by 1680 gives
\begin{equation}
\sum_{(A,B)\in E}p(|A|,|B|)
\le
1680-42N
=
42(40-N).
\label{eq:budget}
\end{equation}
This is the common counting budget for $35\le N\le40$.

If $e$ counts lower sets of sizes $2,5$, each of its two incident edges
has penalty at least 14. The same assertion holds if $e$ instead counts
upper sets of sizes $3,6$. For either family separately,
\begin{equation}
28e\le1680-42N
\quad\Longleftrightarrow\quad
2e+3N\le120.
\label{eq:extreme}
\end{equation}

\section{The analytic exclusion of order 40}
\label{sec:analytic}

The nonnegative penalties in~\eqref{eq:budget} already give $N\le40$.
We now exclude the equality case. This section uses only the chain
identity and a graph argument; none of the finite refutations is needed.

\subsection{A seven-vertex graph obstruction}

\begin{lemma}
\label{lem:graph-seven}
There is no triangle-free graph $G$ on seven vertices such that
$|N_G(u)\cup N_G(v)|=2$ for every nonadjacent pair $u,v$.
\end{lemma}

\begin{proof}
First, every vertex has degree at most two. Indeed, if a vertex $u$
had degree at least three, any nonneighbor $v$ would satisfy
$|N(u)\cup N(v)|\ge3$. Thus $u$ would have to be universal.
Triangle-freeness would then make the other six vertices independent,
and two of them would have a common singleton neighborhood, again
contradicting the required value two.

Next, there are no isolated vertices. If $z$ were isolated, the condition
applied to $z$ and any other vertex $v$ would give $\deg(v)=2$.
The other six vertices would form a triangle-free 2-regular graph,
necessarily a 6-cycle. Two vertices at distance two in that cycle have
a neighborhood union of size three, a contradiction.

It remains to show that $G$ is connected. Suppose it were disconnected.
If one component contained a degree-two vertex $u$, choose a vertex $v$
from another component. There are no isolated vertices, so
$\deg(v)\ge1$. The neighborhoods of $u$ and $v$ lie in different
components, and hence
\[
|N(u)\cup N(v)|
=
\deg(u)+\deg(v)
\ge3,
\]
a contradiction. Therefore every vertex in a disconnected $G$ would
have degree one. Such a graph is a disjoint union of edges and cannot
have seven vertices. Thus $G$ is connected.

A connected graph on seven vertices with maximum degree two is a path
$P_7$ or a cycle $C_7$. Each contains two degree-two vertices at distance
two whose neighborhood union has size three. This final contradiction
proves the lemma.
\end{proof}

\subsection{Excluding equality in the chain bound}

\begin{proposition}[Analytic length bound]
\label{prop:forty}
A saturated inclusion cycle on eight coordinates has at most 39 lower
vertices.
\end{proposition}

\begin{proof}
For $N<25$ there is nothing to prove. Otherwise, the rank bounds and
nonnegative penalties in~\eqref{eq:budget} give $N\le40$.

Suppose $N=40$. The chain budget is then zero. Its zero-penalty edge
types are $(3,4)$ and $(4,5)$, so connectedness forces one homogeneous
type. After complement duality, assume that the lower family
$\mathcal A$ consists of 40 four-sets and the upper family
$\mathcal B$ of 40 five-sets.

The chain identity now has value $40/40=1$.
Every maximal chain must therefore meet $\mathcal A\cup\mathcal B$.
In particular, a four-set outside $\mathcal A$ cannot be contained
in a five-set outside $\mathcal B$: a maximal chain through that
pair would meet neither selected family.

Let $\mathcal F$ be the family of triples complementary to the
five-sets outside $\mathcal B$. Let $\mathcal H$ be the family
of four-sets complementary to members of $\mathcal A$.
If $Q\in\mathcal H$, its complement is a selected four-set,
which has exactly two selected five-set supersets out of its four
supersets. Thus $Q$ contains exactly two members of $\mathcal F$.
If $Q\notin\mathcal H$, all four supersets of its complement are
selected five-sets, so it contains no member of $\mathcal F$.
Every four-set consequently contains either zero or two members
of $\mathcal F$, and belongs to $\mathcal H$ exactly in the latter case.

Also, every triple $T\notin\mathcal F$ is contained in exactly
two members of $\mathcal H$: its complementary five-set belongs
to $\mathcal B$ and has exactly two selected lower neighbors.

Fix a coordinate $v$ and form a graph on the other seven coordinates,
with edge $ab$ when $\{v,a,b\}\in\mathcal F$.
This graph is triangle-free, since a triangle would give a four-set
containing at least three members of $\mathcal F$.

For a nonedge $ab$, the triple $\{v,a,b\}$ is not in $\mathcal F$.
For each remaining coordinate $x$, the four-set $\{v,a,b,x\}$
belongs to $\mathcal H$ precisely when $ax$ or $bx$ is an edge.
Indeed, if neither is an edge, its only possible member of
$\mathcal F$ is $\{a,b,x\}$, and the zero-or-two rule forces
zero. If at least one is an edge, the same rule forces two.

Since $\{v,a,b\}$ has exactly two extensions in $\mathcal H$,
the graph satisfies
\[
|N(a)\cup N(b)|=2
\]
for every nonedge. This contradicts the lemma.
\end{proof}

\section{Reduction to fifteen finite exclusions}
\label{sec:reductions}

\subsection{Selection formulas and valid proper-cycle clauses}

The remaining exclusions use selection formulas with
$35\le|\A|,|\B|\le40$. Each selected vertex has degree two, so the
two selected cardinalities are equal. Write their common value as $N$.
Keeping $N=40$ is a harmless relaxation of the analytic length bound
and preserves the saved CNFs.

Each candidate set has a Boolean variable recording whether it is
selected. A conjunctive normal form (CNF) formula is a conjunction of
clauses, each a disjunction of literals. A checked refutation proves
that no assignment satisfies all clauses. We use reverse-unit-propagation
(RUP) refutations, whose checking rule is given in Section S4 of the
\emph{Supplementary Material}.

The formulas impose the antichain and saturation conditions, the
stage-specific restrictions proved below, and valid clauses excluding
short cycle components. They need not impose full connectedness.
Section S1 of the \emph{Supplementary Material} constructs twenty-one
proper-cycle representatives; Section S2 specifies all encodings.
The next lemma justifies their use without any prior solver verdict.

\begin{lemma}[Proper-cycle clause]
\label{lem:subtour}
Let $z_v$ indicate selection of a set-vertex. If $C$ is an induced
alternating cycle with fewer than 70 vertices in the strict-inclusion
graph, every induced $2N$-cycle with $N\ge35$ satisfies
\begin{equation}
\bigvee_{v\in C}\neg z_v.
\label{eq:subtour}
\end{equation}
The same holds for every coordinate-permuted copy of $C$.
\end{lemma}

\begin{proof}
If all vertices of $C$ were selected, each would already have two
selected neighbors in $C$. The degree-two constraints would preclude any
edge to the remaining selected vertices. Thus $C$ would be a proper
connected component, contradicting connectedness of the target cycle.
Coordinate permutations preserve inclusion.
\end{proof}

The formula names distinguish local branches, not polygon orders.
The following observation links each finite refutation to the covering
problem.

\begin{lemma}[From a cycle to a formula]
\label{lem:certificate-interface}
Fix a normalized local branch. A target cycle satisfying its established
rank restrictions, orientation, and counting inequalities extends to a
satisfying assignment of the corresponding selection formula.
Consequently, a checked refutation of that formula excludes the branch.
\end{lemma}

\begin{proof}
Select exactly the target vertices. Their antichain, degree, and
saturation identities satisfy the base clauses; normalization satisfies
the fixed-vertex clauses. Pruning uses only these conditions and prior
rank exclusions, so it removes no target vertex. The hypotheses give
the counting and orientation constraints, and Lemma~\ref{lem:subtour}
gives the proper-cycle clauses. The encodings in Section S2 of the
\emph{Supplementary Material} then provide the auxiliary values.
A refutation excludes this satisfying assignment.
\end{proof}

\subsection{Stage H: excluding homogeneous cycles}
\label{sec:homogeneous}

For a homogeneous cycle of type $(2,3)$ or $(5,6)$, the extreme count
in~\eqref{eq:extreme} is $e=N$, giving $5N\le120$, impossible here.
Together with the rank restrictions, this leaves the homogeneous types
$(3,4)$ and $(4,5)$. Complement duality reduces them to $(4,5)$.

\begin{lemma}[Five homogeneous branches]
\label{lem:hom-branches}
Every homogeneous cycle of type $(4,5)$ with $N>5$ lower vertices
can be relabeled so that its first three lower sets are
\begin{equation}
A_0=\{0,1,2,3\},
\quad
A_1=\{0,1,2,4\},
\quad
A_2=\{0,1,4,5\}.
\label{eq:hom-prefix}
\end{equation}
Its fourth set is $(A_2\setminus\{r\})\cup\{s\}$, where the pair
$(r,s)$ can be chosen from
\begin{equation}
(4,3),\ (1,3),\ (5,6),\ (4,6),\ (1,6).
\label{eq:hom-five}
\end{equation}
\end{lemma}

\begin{proof}
The intersection of three consecutive lower sets has size two or
three. If every such intersection had size three, start with the
three-element intersection of two consecutive lower sets. Each
successive triple condition preserves this same intersection. All
lower sets would contain a common three-element core, giving at most
$8-3=5$ distinct lower sets, a contradiction.

Choose a triple whose intersection has size two. Normalize the first
two sets to $\{0,1,2,3\}$ and $\{0,1,2,4\}$. The third deletes one
of their common elements and adds an element outside their union.
Adding the previously deleted element 3 instead would repeat the
upper union, which is impossible. A coordinate permutation therefore
gives~\eqref{eq:hom-prefix}.

Write the fourth set as
\[
A_3=(A_2\setminus\{r\})\cup\{s\},
\qquad
r\in\{0,1,4,5\},
\quad
s\in\{2,3,6,7\}.
\]
The earlier upper set $B_1=\{0,1,2,4,5\}$ would contain $A_3$ if
$s=2$, so this value is forbidden. For $s\in\{3,6,7\}$, the new
upper union $A_2\cup\{s\}$ omits coordinate 2 and therefore contains
neither $A_0$ nor $A_1$. The only further obstruction is containment
in $B_0=\{0,1,2,3,4\}$: it occurs precisely when $(r,s)=(5,3)$.
Thus the complete condition is
\begin{equation}
s\in\{3,6,7\},
\qquad
(r,s)\ne(5,3).
\label{eq:fourth-set-rule}
\end{equation}

The prefix stabilizer independently swaps $0,1$ and $6,7$.
For $s=3$, it distinguishes deletion of a core element from deletion
of 4; deletion of 5 is forbidden. For $s\in\{6,7\}$, it distinguishes
deletion of a core element, of 4, and of 5. This gives exactly five
orbits, represented by~\eqref{eq:hom-five}.
\end{proof}

\begin{lemma}[Computational lemma H]
\label{lem:H}
For $35\le N\le40$, no homogeneous saturated inclusion cycle of type
$(4,5)$ exists on eight coordinates.
\end{lemma}

\begin{proof}
For each of the five branches in~\eqref{eq:hom-five}, fix the four lower
sets and their three consecutive upper unions. Use all four-sets as
lower candidates and all five-sets as upper candidates, with selection
sizes between 35 and 40, and impose the base clauses and the proper-cycle
clauses described above.

These are the formulas \code{H43}, \code{H57}, \code{H83},
\code{H99}, and \code{H113}. Each has a checked RUP refutation;
their sizes and proof lengths are recorded in Table S2 of the
\emph{Supplementary Material}.
Lemma~\ref{lem:hom-branches} covers all possibilities.
\end{proof}

Complement duality also excludes the homogeneous type $(3,4)$.

\subsection{Stage G: restricting the possible rank pairs}
\label{sec:rank-exclusions}

We next exclude three rank pairs and their duals. These exclusions
also sharpen the counting budget for the subsequent stages.

\begin{lemma}[Computational lemma G]
\label{lem:G}
A saturated inclusion cycle with $35\le N\le40$ lower vertices has no
edge of type $(2,3)$, $(2,5)$, or $(2,6)$.
\end{lemma}

\begin{proof}
For each chosen type $(2,s)$, normalize its edge to
\[
A=\{0,1\}
\subset
B=\{0,\ldots,s-1\}.
\]
The other lower neighbor $C$ of $B$ satisfies $A\cup C=B$.
Thus
\[
C=(B\setminus A)\cup X
\]
with $X$ a proper subset of $A$.
The lower rank bound and the swap of 0 and 1 give these representatives:
\begin{equation}
\begin{array}{c|c|c}
(2,s)&B\setminus A&X\\
\hline
(2,3)&\{2\}&\{0\}\\
(2,5)&\{2,3,4\}&\varnothing\ \text{or}\ \{0\}\\
(2,6)&\{2,3,4,5\}&\varnothing\ \text{or}\ \{0\}
\end{array}
\label{eq:G-anchors}
\end{equation}

Use the full rank ranges~\eqref{eq:rank-range}, fix the chosen edge,
require at least one of these representative neighbors, and impose
the base clauses,~\eqref{eq:extreme}, and the proper-cycle clauses.
The resulting formulas \code{G23}, \code{G25}, and \code{G26} each have
a checked RUP refutation. The normalization uses only the set ranks,
so it is valid uniformly in $N$.
\end{proof}

\begin{corollary}
\label{cor:five-types}
Every remaining saturated inclusion cycle has lower ranks in
$\{2,3,4\}$ and upper ranks in $\{4,5,6\}$. Its possible edge types
are exactly among
\begin{equation}
(2,4),\ (3,4),\ (3,5),\ (4,5),\ (4,6).
\label{eq:five-types}
\end{equation}
\end{corollary}

\begin{proof}
Complement duality also excludes $(5,6)$ and $(3,6)$.
A lower set of size five could only have upper neighbors of size six,
and an upper set of size three could only have lower neighbors of size
two. Both are impossible by these exclusions. The remaining strict
rank pairs give~\eqref{eq:five-types}.
\end{proof}

Let $x,a,u$ count lower sets of sizes $2,3,4$, and let $b,z$ count upper
sets of sizes $4,6$, respectively. Then $N=x+a+u$.
Write $t$ for the number of $(3,5)$ edges.
All $2x$ edges incident with lower two-sets end at upper four-sets.
Consequently the number of $(3,4)$ edges is $2b-2x$, and counting at
lower three-sets gives
\begin{equation}
t=2a-(2b-2x)=2(a-b+x).
\label{eq:t-identity}
\end{equation}

The only nonzero penalties among~\eqref{eq:five-types} are 17 on
$(2,4),(4,6)$ and 6 on $(3,5)$. Therefore
\begin{align}
34(x+z)+6t
&\le1680-42N,
\label{eq:refined-budget}\\
23x+17z+6a-6b+21N
&\le840,
\label{eq:linearN}\\
44x+27a+21u+17z-6b
&\le840.
\label{eq:constant}
\end{align}
The last identity uses $N=x+a+u$. It is the constant-right-hand-side
weighted constraint used in the uniform formulas. Since $t\ge0$ and
$N\ge35$, the first inequality also gives
\begin{equation}
x+z\le6.
\label{eq:six}
\end{equation}

\subsection{Stage T: removing lower two-sets and upper six-sets}
\label{sec:twos}

After G, the extreme ranks are lower two-sets and upper six-sets.
Duality permits us to anchor at a lower two-set. One of its four local
branches is excluded by the following shadow argument; the other three
are excluded by checked refutations.

\begin{lemma}[Three disjoint pairs]
\label{lem:three-pairs}
Suppose an inclusion cycle in $2^{\Om_8}$ has all upper ranks at least
four. If its lower family contains three pairwise disjoint two-element
sets, then its upper family has at most 34 members.
\end{lemma}

\begin{proof}
Let the three pairs be $P_0,P_1,P_2$. Split the upper family into
$\mathcal X$, whose members contain at least one $P_i$, and
$\mathcal Y$, whose members contain none. Each $P_i$ has exactly two
upper neighbors, so there are six pair--upper incidences and
$|\mathcal X|\le6$.

A set avoiding all three pairs chooses at most one element from each
pair and can also use either of the two remaining coordinates. Its
size is at most five. The number of such four-sets is
\begin{equation}
2^3\binom21+\binom32\,2^2
=
16+12
=
28.
\label{eq:independent-four-count}
\end{equation}
The first term uses all three pairs and one remaining coordinate;
the second uses two pairs and both remaining coordinates.

Write $\mathcal Y_j$ for the members of $\mathcal Y$ of size $j$.
Let $\Delta$ be the family of four-subsets of members of
$\mathcal Y_5$. Every member of $\mathcal Y_5$ has five such subsets.
Each such four-set extends to at most two such five-sets.
If it uses all three pairs, only the second free coordinate can be
added. If it uses two pairs, either element of the missing pair can
be added.
Counting containments gives
\[
5|\mathcal Y_5|\le2|\Delta|.
\]

Since the upper family is an antichain,
$\mathcal Y_4\cap\Delta=\varnothing$. Therefore
\begin{equation}
|\mathcal Y|
\le
28-|\Delta|+|\mathcal Y_5|
\le
28-\frac32|\mathcal Y_5|
\le28.
\label{eq:three-pair-shadow}
\end{equation}
It follows that
\[
|\B|=|\mathcal X|+|\mathcal Y|\le6+28=34.
\]
\end{proof}

\begin{lemma}[Four local branches]
\label{lem:two-branches}
If $x+z>0$, complement duality permits $x\ge z$.
After relabeling, a five-vertex subpath has the form
\begin{equation}
P_1\cup X,\quad
P_0\cup P_1,\quad
P_0,\quad
P_0\cup P_2,\quad
P_2\cup Y,
\label{eq:two-path}
\end{equation}
where
\[
P_0=\{0,1\},
\qquad
P_1=\{2,3\},
\qquad
P_2=\{4,5\},
\]
and
\begin{equation}
\begin{array}{c|cc}
\text{branch}&X&Y\\
\hline
\mathrm T_0&\varnothing&\varnothing\\
\mathrm T_1&\varnothing&\{0\}\\
\mathrm T_2&\{0\}&\{0\}\\
\mathrm T_3&\{0\}&\{1\}
\end{array}
\label{eq:two-four}
\end{equation}
\end{lemma}

\begin{proof}
Duality exchanges $x$ and $z$, so choose $x\ge z$. Then $x\ge1$.
Fix a two-element lower set $P_0$. Both upper neighbors have size
four, and their intersection is precisely $P_0$. Their two-element
parts outside $P_0$ are therefore disjoint; call them $P_1,P_2$.
By the union identity, the other lower neighbor on each side is
$P_1\cup X$ or $P_2\cup Y$, where $X,Y$ are empty or singleton
subsets of $P_0$.

Up to simultaneously exchanging 0 and 1 and reversing the path,
the complete possibilities are: both empty; exactly one empty;
two equal singletons; or two distinct singletons. These give the
four rows of~\eqref{eq:two-four}.
\end{proof}

\begin{lemma}[Computational lemma T]
\label{lem:T}
Every remaining saturated inclusion cycle has $x=z=0$.
\end{lemma}

\begin{proof}
If $x+z>0$, choose complement orientation $x\ge z$ and apply
Lemma~\ref{lem:two-branches}. Branch $\mathrm T_0$ contains three
pairwise disjoint lower pairs, contradicting Lemma~\ref{lem:three-pairs}
because $N\ge35$.

For each other branch, fix its five-vertex path and use the remaining
rank pairs, the base and proper-cycle clauses,
\eqref{eq:constant}, \eqref{eq:six}, and $x\ge z$.
The three formulas \code{T1}, \code{T2}, and \code{T3} have checked RUP
refutations. Thus none of the four branches occurs.
\end{proof}

\subsection{Stage M: proving that mixed edges form a matching}
\label{sec:middle}

Only lower ranks $3,4$ and upper ranks $4,5$ remain.
Call an edge of type $(3,5)$ mixed. If there were no mixed edge,
connectedness would force a homogeneous type $(3,4)$ or $(4,5)$,
contrary to Lemma~\ref{lem:H} and duality. Hence $t>0$, and
\begin{equation}
t=2(a-b),
\qquad
6t\le1680-42N.
\label{eq:middle-budget}
\end{equation}

All preceding exclusions are preserved by complement duality.
Under this duality, $(a,b)$ becomes $(N-b,N-a)$.
Choose the orientation
\begin{equation}
a+b\le N,
\qquad\text{equivalently}\qquad
b\le u.
\label{eq:middle-orientation}
\end{equation}
This choice precedes the local normalization, which uses only coordinate
permutations. The whole cycle cannot be mixed: $t=2N$ would imply
$54N\le1680$, contrary to $N\ge35$.

\begin{lemma}[Three boundary branches]
\label{lem:middle-branches}
In a cycle satisfying~\eqref{eq:middle-orientation}, a mixed edge can
be chosen and normalized as
\[
A=\{0,1,2\}
\subset
B=\{0,1,2,3,4\}.
\]
Let $C$ be the other lower neighbor of $B$ and $D$ the other upper
neighbor of $A$. Write
\[
C=\{3,4\}\cup X,
\qquad
D=A\cup Y.
\]
Then the three possible branches have representatives
\begin{equation}
\begin{array}{c|cc}
\text{branch}&X&Y\\
\hline
\mathrm M_0&\{0\}&\{5\}\\
\mathrm M_1&\{0,1\}&\{5\}\\
\mathrm M_2&\{0,1\}&\{5,6\}
\end{array}
\label{eq:middle-three}
\end{equation}
\end{lemma}

\begin{proof}
By~\eqref{eq:middle-budget}, mixed edges do not occupy the entire
$2N$-cycle. Choose a mixed edge at a boundary of a maximal run of mixed
edges. At least one of its adjacent edges is not mixed.

Normalize $A=\{0,1,2\}$ and $B=\{0,1,2,3,4\}$. The union identity
requires that $C$ contain $\{3,4\}$ and either one or two elements of
$A$. The intersection identity requires that $D$ be $A$ together with
either one or two elements of $\{5,6,7\}$. The core coordinates and
the outside coordinates may be permuted independently.

If $|C|=3$ and $|D|=5$, both adjacent edges are mixed, contrary to the
boundary choice. For each remaining pair of sizes, the independent
permutations of $A$ and $\{5,6,7\}$ act transitively on the choices.
The three size pairs therefore give exactly the representatives
in~\eqref{eq:middle-three}.
\end{proof}

The budget and orientation yield three further necessary bounds:
\begin{equation}
1\le a-b\le17,
\qquad
a\le26,
\qquad
b\le19.
\label{eq:bounds}
\end{equation}
Indeed, $t=2(a-b)>0$ and
\[
6(a-b)+21N\le840
\]
give the first bounds.
Since $b\le N-a$, the same inequality gives
\[
12a+15N\le840,
\]
and $N\ge35$ gives $a\le26$.
Finally, $a\ge b+1$ and $a+b\le N\le40$ imply $b\le19$.
These consequences are valid before imposing any matching condition.

\begin{lemma}[Computational lemma M]
\label{lem:M}
Neither boundary branch $\mathrm M_0$ nor $\mathrm M_2$ occurs in a
remaining cycle satisfying~\eqref{eq:middle-orientation}.
\end{lemma}

\begin{proof}
Use lower ranks $3,4$ and upper ranks $4,5$, and fix the four vertices
of the indicated normalized branch. Impose the base and proper-cycle
clauses,~\eqref{eq:constant} with $x=z=0$,
\eqref{eq:middle-orientation}, and~\eqref{eq:bounds}.
The formulas \code{M0} and \code{M2} have checked RUP refutations.
\end{proof}

\begin{corollary}
\label{cor:matching}
The mixed edges of a remaining cycle form a nonempty matching.
\end{corollary}

\begin{proof}
If a maximal run of mixed edges had length at least two, consider an
endpoint mixed edge of that run. If its nonmixed neighbor is incident
with its lower endpoint, then $|D|=4$ and $|C|=3$, giving
$\mathrm M_0$. If the nonmixed neighbor is incident with its upper
endpoint, then $|C|=4$ and $|D|=5$, giving $\mathrm M_2$.
Both are impossible by Lemma~\ref{lem:M}. A run occupying the entire
cycle has already been excluded by~\eqref{eq:middle-budget}. Thus
every run has length one. Nonemptiness was proved above.
\end{proof}

In particular, it is now valid to require that a selected rank-three
lower set have at most one selected rank-five upper neighbor, and
that a selected rank-five upper set have at most one selected
rank-three lower neighbor. These additional constraints depend on
\emph{both} refutations in Lemma~\ref{lem:M}; they are not assumptions
used in proving that lemma.

\subsection{Stage F: excluding the final two branches}
\label{sec:final}

The matching property has now been proved. It reduces the neighborhood
of any mixed edge to the following two possibilities.

\begin{lemma}[Two final branches]
\label{lem:final-branches}
A remaining cycle can be normalized to contain lower sets
\[
A=\{0,1,2\},
\qquad
C=\{0,1,3,4\}
\]
and upper sets
\[
B=\{0,1,2,3,4\},
\qquad
D=\{0,1,2,5\}.
\]
The other lower neighbor $E$ of $D$ can be chosen as
$\{0,1,5\}$ or $\{0,2,5\}$.
\end{lemma}

\begin{proof}
Choose any mixed edge. By Corollary~\ref{cor:matching}, its adjacent
edges are both nonmixed. It is therefore in branch $\mathrm M_1$,
so normalize its four sets as in~\eqref{eq:middle-three}.

The other lower neighbor $E$ of $D$ has size three and satisfies
$A\cup E=D$. Its possibilities are exactly
\[
\{0,1,5\},
\qquad
\{0,2,5\},
\qquad
\{1,2,5\}.
\]
The swap $0\leftrightarrow1$ fixes $A,B,C,D$, fixes the first
possibility, and exchanges the last two. Thus two representatives
suffice.
\end{proof}

\begin{lemma}[Computational lemma F]
\label{lem:F}
Neither final branch in Lemma~\ref{lem:final-branches} occurs.
\end{lemma}

\begin{proof}
For $E=\{0,1,5\}$ or $E=\{0,2,5\}$, fix the five normalized vertices.
Use the M-stage constraints and additionally require the mixed edges
to form a matching, as justified by Corollary~\ref{cor:matching}.
For a vertex $v$ and two of its mixed neighbors $r,s$, this last
condition is the clause
\[
\neg z_v\lor\neg z_r\lor\neg z_s.
\]
The two formulas \code{F35} and \code{F37} have checked RUP refutations.
The suffixes are the decimal masks of $E$, not matrix orders.
\end{proof}

\section{Completion of the main proof}
\label{sec:main-proof}

\begin{proof}[Proof of Theorem~\ref{thm:uniform}]
Suppose $35\le n\le40$ and a covering by at most eight rectangles exists.
We apply the reductions in the order described in Section~\ref{sec:strategy}.

\medskip
\noindent\textbf{Step 1: normalization and the analytic bound.}
Pad the covering with empty rectangles and apply Lemma~\ref{lem:normal}.
The resulting saturated inclusion cycle has $N=n$ lower vertices.
Proposition~\ref{prop:ranks} gives its rank ranges, and
\eqref{eq:budget} gives its penalty budget. Proposition~\ref{prop:forty}
gives $N\le39$, leaving only $35\le N\le39$.

\medskip
\noindent\textbf{Step 2: the H, G, and T exclusions.}
Lemma~\ref{lem:H} and duality exclude both homogeneous types.
Lemma~\ref{lem:G} and duality leave the five rank pairs
in~\eqref{eq:five-types}. Lemma~\ref{lem:three-pairs} and the three
T refutations then give $x=z=0$ by Lemma~\ref{lem:T}.
Thus the lower ranks are $3,4$ and the upper ranks are $4,5$.

\medskip
\noindent\textbf{Step 3: the M exclusions.}
Since the cycle is not homogeneous, it contains a mixed edge.
Choose the dual orientation~\eqref{eq:middle-orientation}.
Lemma~\ref{lem:M} excludes the two boundary branches that would allow
adjacent mixed edges. Corollary~\ref{cor:matching} therefore gives a
nonempty matching of mixed edges.

\medskip
\noindent\textbf{Step 4: the final contradiction.}
Lemma~\ref{lem:final-branches} reduces a neighborhood of a mixed edge
to one of two F branches. Both are excluded by Lemma~\ref{lem:F}.
Every finite exclusion is applicable by
Lemma~\ref{lem:certificate-interface}, using only the conditions
established before that stage.

The assumed covering cannot exist, so $\RCB(M_n)\ge9$.
Proposition~\ref{prop:known-upper} supplies the published upper bound
$\RCB(M_n)\le9$. Since rectangle covering number depends only on the
support, the same equality holds for $S_P$.
\end{proof}

\section{Certificates and supplementary material}
\label{sec:certificates}

All fifteen CNFs were reconstructed byte for byte. A RUP implementation
separate from the solver checked their refutations through an empty
clause, verifying 7,023,863 additions in total.
Lemma~\ref{lem:certificate-interface} and the local classifications give
the mathematical reduction; the refutations give the finite exclusions.
The analytic bound in Proposition~\ref{prop:forty} uses none of these
refutations.

The accompanying \emph{Supplementary Material} provides the proper-cycle
constructions and complete parameters in Section S1, the exact encodings
in Section S2, the certificate inventory in Section S3, checking methods
and commands in Section S4, and the published upper-bound witnesses in
Section S5. It is supplied together with this paper and the certificate
archive. The checks reported here do not assert a third-party audit or
a proof-assistant formalization.

The certificate archives are available at
\url{https://doi.org/10.5281/zenodo.23837747}.

\section*{Declaration of generative-AI assistance}
\addcontentsline{toc}{section}{Declaration of generative-AI assistance}

An OpenAI language model accessed through Codex assisted with the 
programs, certificate workflow, and writing.
Its assertions were not treated as certificates. The finite exclusions
were checked by the supplied RUP implementation.

\clearpage
\section*{Supplementary Material}
\addcontentsline{toc}{section}{Supplementary Material}
\setcounter{section}{0}
\setcounter{subsection}{0}
\setcounter{table}{0}
\setcounter{equation}{0}
\renewcommand{\thesection}{S\arabic{section}}
\renewcommand{\thetable}{S\arabic{table}}
\numberwithin{equation}{section}
\renewcommand{\theHsection}{supp.\arabic{section}}
\renewcommand{\theHtable}{supp.\arabic{table}}
\renewcommand{\theHequation}{supp.\arabic{section}.\arabic{equation}}
\providecommand{\theHtheorem}{}
\renewcommand{\theHtheorem}{supp.\arabic{section}.\arabic{theorem}}

This supplement accompanies the main paper of the same title and is
provided as part of its proof record. It collects the complete
proper-cycle parameters, encoding details, certificate inventory,
checking procedures, and published upper-bound witnesses.
The main paper retains the key mathematical lemmas and local
classification proofs. Neither source requires the other document's
auxiliary files for compilation.

We use $M_n$ for the zero--one matrix whose zeros are $(j,j)$ and
$(j+1,j)$ modulo $n$. A saturated inclusion cycle consists of two
antichains $\A=\{A_i\}$ and $\B=\{B_j\}$ satisfying
\[
 B_j=A_j\cup A_{j+1},\qquad
 A_i=B_{i-1}\cap B_i,\qquad
 A_i\subseteq B_j\Longleftrightarrow i\in\{j,j+1\}.
\]
The saved formulas use $\Omega_8=\{0,\ldots,7\}$ and allow
$N=|\A|=|\B|$ with $35\le N\le40$. The analytic length bound in
the main paper excludes $N=40$ without a refutation. Keeping this
value in the formulas is a harmless relaxation and preserves the
originally checked CNFs.

\section{Proper-cycle constructions and complete parameters}
\label{app:supp:cycles}

The proper-cycle clauses belong to the initial CNFs. They must
therefore be justified independently of any solver verdict. The following three constructions establish the validity of
every supplied representative without testing its incidence matrix.

\subsection{Why proper-cycle clauses are valid}

\begin{lemma}[Proper-cycle clause]\label{app:lem:subtour}
Let $z_v$ indicate selection of a set-vertex. If $C$ is an induced
alternating cycle with fewer than 70 vertices in the strict-inclusion
graph, every induced $2N$-cycle with $N\ge35$ satisfies
\begin{equation}
 \bigvee_{v\in C}\neg z_v.
 \label{app:eq:subtour}
\end{equation}
The same holds for every coordinate-permuted copy of $C$.
\end{lemma}

\begin{proof}
If all vertices of $C$ were selected, each would already have two
selected neighbors in $C$. The degree-two constraints would preclude any
edge to the remaining selected vertices. Thus $C$ would be a proper
connected component, contradicting connectedness of the target cycle.
Coordinate permutations preserve inclusion.
\end{proof}

\subsection{General cycle constructions}

\begin{lemma}[Cyclic intervals]\label{app:lem:interval-cycles}
Let $u_0,\ldots,u_{m-1}$ be distinct coordinates, let $T$ be disjoint
from them, and let $1\le s\le m-2$. With indices modulo $m$, define
\[
 L_i=T\cup\{u_i,\ldots,u_{i+s-1}\},\qquad
 U_i=T\cup\{u_i,\ldots,u_{i+s}\}.
\]
Then $(L_i,U_i)$ is a saturated inclusion cycle. Denote this
construction by $I_{m,s}(T;u_0,\ldots,u_{m-1})$.
\end{lemma}

\begin{proof}
The union and intersection identities follow from consecutive
intervals and $s\le m-2$. To check inclusion, fix $j$ and regard
$\{u_j,\ldots,u_{j+s}\}$ as a proper arc of the cyclic order.
A consecutive $s$-element interval contained in this arc can start
only at $u_j$ or $u_{j+1}$. A later start would meet $u_{j+s+1}$,
which is outside the arc; a start outside the arc is immediately
excluded. Hence $L_i\subseteq U_j$ exactly when $i=j$ or $j+1$.
\end{proof}

\begin{lemma}[Subdivision using fresh coordinates]\label{app:lem:subdivision}
Suppose consecutive lower sets in a saturated inclusion cycle are
$P\cup\{r\}$ and $P\cup\{s\}$, with $r\ne s$ outside $P$.
For fresh coordinates $z_1,\ldots,z_q$ absent from every old set,
replace this adjacency by the lower-set path
\[
 P\cup\{r\},\ P\cup\{z_1\},\ldots,
 P\cup\{z_q\},\ P\cup\{s\},
\]
and use consecutive unions as the new upper sets. The result is
again a saturated inclusion cycle.
\end{lemma}

\begin{proof}
Each new lower set contains a fresh coordinate and therefore lies
in no old upper set. Each new upper set contains either two fresh
coordinates or one fresh coordinate and one of $r,s$. An old lower
set contained in it would consequently lie in $P$, $P\cup\{r\}$,
or $P\cup\{s\}$. The antichain property allows only the appropriate
endpoint in the latter two cases and no old lower set in the first.
Among new sets, the fresh coordinate specifies exactly the two
incident upper sets. Union and intersection identities hold along
the path, and at its endpoints they follow because old sets omit
the fresh coordinates. All other incidences are unchanged.
\end{proof}

\subsection{The H representatives}

For the only subdivision pattern needed here, take eight distinct
coordinates in the order $(t,a,b,c,d,e,x,y)$, put $T=\{t\}$, and
use the cyclic sequence
\begin{equation}
\begin{split}
 T\cup\{a,x,y\},\quad T\cup\{b,x,y\},\quad
 T\cup\{c,x,y\},\quad T\cup\{c,d,x\},\\
 T\cup\{c,d,e\},\quad T\cup\{d,e,y\},\quad
 T\cup\{e,x,y\}.
\end{split}
\label{app:eq:expanded-cycle}
\end{equation}
Call it $E(t,a,b,c,d,e,x,y)$. Start with
$I_{5,3}(\{t\};y,x,c,d,e)$ and subdivide the adjacency between
$T\cup\{e,x,y\}$ and $T\cup\{c,x,y\}$ using the fresh coordinates
$a,b$. Lemma~\ref{app:lem:subdivision} proves~\eqref{app:eq:expanded-cycle}.

Table~\ref{app:tab:cycles} gives exact construction parameters for all
seventeen H representatives. In an $I$ row the third column is $T$
and the fourth is the ordered cyclic list of $u$ coordinates. In an
$E$ row the fourth column is $(t,a,b,c,d,e,x,y)$.
Substitution yields the same selected lower and upper sets as the
saved representatives, up to rotation or reversal of their listing.
The full coordinate orbits, and therefore their proper-cycle clauses,
are unchanged.

\begin{longtable}{@{}lllp{73mm}@{}}
\caption{Analytic descriptions of all H representatives. Coordinates,
unlike decimal masks, are the integers 0 through 7.}
\label{app:tab:cycles}\\
\toprule
ID & Construction & Core $T$ & Ordered coordinates\\
\midrule
\endfirsthead
\toprule ID & Construction & Core $T$ & Ordered coordinates\\\midrule
\endhead
\bottomrule
\endfoot
H43.1 & $I_{5,1}$ & $\{2,3,5\}$ & $0,4,1,6,7$\\
H43.2 & $E$ & -- & $6,4,3,7,1,2,5,0$\\
H43.3 & $I_{5,3}$ & $\{1\}$ & $0,4,7,3,6$\\
H57.1 & $E$ & -- & $0,1,3,7,4,2,5,6$\\
H57.2 & $I_{7,3}$ & $\{7\}$ & $0,4,2,3,1,6,5$\\
H57.3 & $I_{5,3}$ & $\{0\}$ & $1,2,3,5,4$\\
H57.4 & $I_{5,1}$ & $\{2,4,6\}$ & $0,5,1,3,7$\\
H83.1 & $I_{5,2}$ & $\{2,4\}$ & $0,6,5,3,7$\\
H83.2 & $I_{5,3}$ & $\{0\}$ & $2,5,3,4,6$\\
H83.3 & $E$ & -- & $0,5,6,7,3,2,4,1$\\
H99.1 & $I_{5,3}$ & $\{0\}$ & $2,5,3,4,6$\\
H99.2 & $I_{7,3}$ & $\{7\}$ & $0,3,2,5,6,1,4$\\
H99.3 & $I_{5,1}$ & $\{2,4,6\}$ & $0,5,1,3,7$\\
H113.1 & $E$ & -- & $7,2,5,6,0,3,4,1$\\
H113.2 & $I_{5,1}$ & $\{1,5,7\}$ & $0,2,6,4,3$\\
H113.3 & $I_{7,3}$ & $\{0\}$ & $1,6,2,5,3,4,7$\\
H113.4 & $I_{5,3}$ & $\{4\}$ & $1,2,6,5,7$\\
\end{longtable}

\subsection{The Q representatives}

\begin{lemma}[Disjoint petals]\label{app:lem:petals}
Let $P_0,P_1,P_2$ be nonempty pairwise disjoint sets, disjoint from
a core $T$. The lower sets $L_i=T\cup P_i$ and upper sets
$U_i=T\cup P_i\cup P_{i+1}$, with indices modulo 3, form a
saturated six-vertex inclusion cycle.
\end{lemma}

\begin{proof}
Disjointness and nonemptiness give $L_i\subseteq U_j$ exactly when
$i\in\{j,j+1\}$. Taking the indicated unions and intersections
recovers $U_i$ and $L_i$ respectively.
\end{proof}

The four Q representatives follow immediately from these petals:
\[
\begin{array}{c|c|ccc}
 &T&P_0&P_1&P_2\\\hline
 Q1&\{0,1\}&\{2\}&\{3\}&\{4\}\\
 Q2&\{0,1,2\}&\{3\}&\{4\}&\{5\}\\
 Q3&\{0,1\}&\{2\}&\{3\}&\{4,5\}\\
 Q4&\{0\}&\{1,2\}&\{3,4\}&\{5,6\}.
\end{array}
\]
Their lower-mask sequences are respectively $(7,11,19)$,
$(15,23,39)$, $(7,11,51)$, and $(7,25,97)$.

All fifteen uniform formulas use all twenty-one representatives, together
with their coordinate-permuted copies that lie in the retained candidate
domain. Each representative has 6, 10, or 14 vertices, so
Lemma~\ref{app:lem:subtour} applies for every $35\le N\le40$.
The generators compare the saved data exactly with the parameters above.
Thus every proper-cycle clause has a direct mathematical justification.

\section{Selection model and exact encodings}\label{app:supp:encoding}

For reference, $x,a,u$ count selected lower sets of sizes $2,3,4$,
and $b,z$ count selected upper sets of sizes $4,6$. After the G
exclusions, $N=x+a+u$. The main paper proves the following constraints,
which we restate to specify the encodings without external equation
references:
\begin{equation}\label{app:eq:extreme}
 2e+3N\le120.
\end{equation}
Here $e$ counts lower sets of sizes $2,5$ or, separately, upper sets of
sizes $3,6$. After the G exclusions, the weighted constraint is
\begin{equation}\label{app:eq:constant}
 44x+27a+21u+17z-6b\le840.
\end{equation}
After the T exclusions, the M and F stages have $x=z=0$ and use the
dual orientation $a+b\le N$, equivalently $b\le u$, as well as
\begin{equation}\label{app:eq:bounds}
 1\le a-b\le17,\qquad a\le26,\qquad b\le19.
\end{equation}
The matching condition is imposed only at the F stage, after the two
M-stage refutations have justified it in the main paper.

\subsection{Selection variables and saturation}
Encode a set $S\subseteq\Omega_8$ by the decimal mask
$\mu(S)=\sum_{t\in S}2^t$. This explains the numerical suffixes of
the formula names. For each retained lower or upper candidate introduce
a selection variable $z_{(\mathrm a,S)}$ or $z_{(\mathrm b,T)}$.
The base clauses impose the following conditions:
\begin{enumerate}
 \item Select between 35 and 40 sets in each family.
 \item Do not select two distinct comparable sets in either family,
       or the same set in both families.
 \item Every selected vertex has exactly two selected permitted
       strict-inclusion neighbors. For a previously excluded rank pair,
       prohibit selection of both endpoints.
 \item The two lower neighbors of a selected upper set have that set
       as their union; the two upper neighbors of a selected lower set
       have that set as their intersection.
\end{enumerate}
The last requirement has the clauses
\begin{align}
 &\neg z_{(\mathrm b,T)}\lor
   \bigvee_{\substack{S\subsetneq T\\t\in S}}z_{(\mathrm a,S)},
   &&t\in T,\label{app:eq:union-clause}\\
 &\neg z_{(\mathrm a,S)}\lor
   \bigvee_{\substack{T\supsetneq S\\t\notin T}}z_{(\mathrm b,T)},
   &&t\notin S.\label{app:eq:intersection-clause}
\end{align}
The disjunctions range over retained permitted neighbors. Together with
the degree constraints, they express exactly the required set identities.

Cardinality constraints use the sequential-counter encoding supplied by
Python-SAT~\cite{app:PySAT}. Every clause of a selected vertex's neighbor
constraint is guarded by the negation of its selection variable.
All auxiliary variables are fresh for their respective encodings.
If the vertex is unselected, every guarded clause is true. If it is
selected, the clauses are precisely an encoding of degree two.
Thus guarding preserves the intended conditional constraint.

Pruning removes only candidates forced false by an antichain conflict,
a forbidden inclusion rank, or the two neighbors already fixed at a
selected vertex. If a forced vertex has only two remaining possible
neighbors, both are forced selected. Each recorded step is checked
independently. No exact-order cardinality clause from an earlier
$n=35$ builder is retained: the uniform CNF is constructed afresh with
the interval cardinalities stated above.

This construction requires only the implication from a target cycle to
a satisfying assignment. A satisfying assignment may describe several
components; it need not give a polygon covering. Such extra assignments
make unsatisfiability harder to establish and do not invalidate an
exclusion.

The G formulas use~\eqref{app:eq:extreme}. The T formulas use
\eqref{app:eq:constant}, $x+z\le6$, and $x\ge z$.
The M and F formulas use~\eqref{app:eq:constant} with $x=z=0$,
$a+b\le N$, and~\eqref{app:eq:bounds}.
Only the F formulas impose the matching condition, after both M
refutations have established the matching property in the main paper.
The absence of full connectivity clauses makes these formulas
relaxations of the target cycle problem; all added proper-cycle clauses
are valid by Lemma~\ref{app:lem:subtour}. Unsatisfiability of these relaxations
therefore suffices.

\subsection{The weighted encoding}
\begin{proposition}\label{app:prop:adder-encoding}
A network of Boolean adders followed by the binary comparison clauses
below is an existential CNF encoding of a nonnegative integer weighted
inequality $\sum_j w_j v_j\le B$.
\end{proposition}

\begin{proof}
A half-adder satisfies $u+v=s+2c$, where
$s=u\mathbin\oplus v$ and $c=u\land v$. A full-adder satisfies
\[
 u+v+w=s+2c,\qquad
 s=u\mathbin\oplus v\mathbin\oplus w,\quad
 c=(u\land v)\lor(u\land w)\lor(v\land w).
\]
The sum output is the parity, and the carry output states whether
at least two inputs are one. These identities prove the arithmetic
of each gate. Clauses defining the displayed Boolean equivalences
give an exact CNF representation of the gate.

Expand each $w_j$ in binary and place its input literal $v_j$ in the
corresponding bit columns. Repeatedly replace two or three column
inputs by a sum bit in that column and a carry bit in the next.
Summing the gate identities with their powers-of-two weights shows
that the resulting binary word $S$ equals $\sum_jw_jv_j$.

Pad $S$ and $B$ with leading zeroes to a common length. For every
position $i$ with $B_i=0$, use the clause
\begin{equation}
 \neg S_i\ \lor\!
 \bigvee_{\substack{j>i\\B_j=1}}\neg S_j
 \ \lor\!
 \bigvee_{\substack{j>i\\B_j=0}} S_j.
 \label{app:eq:binary-comparator}
\end{equation}
It forbids $S_i=1$ when all more significant bits match $B$.
If $S>B$, the highest differing bit violates exactly such a clause.
If $S\le B$, every such clause is satisfied, either at position $i$
or by a higher differing bit. Thus these clauses express $S\le B$.
The gate equations determine their auxiliary values from the
inputs, proving existential equivalence.
\end{proof}

In~\eqref{app:eq:constant}, let $m_b$ be the number of retained upper
four-set candidates. Replace each $-6z_{(\mathrm b,B)}$ by
$6(1-z_{(\mathrm b,B)})$ and increase the right-hand side to
$840+6m_b$. This identity reduces the signed inequality to the
nonnegative weighted encoding in Proposition~\ref{app:prop:adder-encoding}.

\section{Certificate inventory and verification record}\label{app:supp:certificates}

Table~\ref{app:tab:certificates} records the fifteen uniform formulas.
Each allows the selected cardinality to vary over the entire interval
$35\le N\le40$. Every formula has been reconstructed byte for byte,
and every refutation has been checked through an empty clause by the
separate RUP implementation supplied with this work.

\begin{table}[htbp]\centering\small
\begin{tabular}{@{}lrrr@{}}\toprule
Case & Variables & Initial clauses & Checked RUP additions\\\midrule
\code{H43} & 5,158 & 38,915 & 316,971\\
\code{H57} & 5,158 & 38,915 & 325,151\\
\code{H83} & 5,158 & 38,915 & 430,293\\
\code{H99} & 5,158 & 38,915 & 206,040\\
\code{H113} & 5,158 & 38,915 & 285,331\\
\code{G23} & 42,272 & 104,461 & 339,780\\
\code{G25} & 46,372 & 115,028 & 501,582\\
\code{G26} & 43,730 & 105,857 & 132,152\\
\code{T1} & 13,402 & 33,498 & 31,624\\
\code{T2} & 18,444 & 47,655 & 270,546\\
\code{T3} & 18,032 & 46,516 & 254,984\\
\code{M0} & 27,355 & 69,151 & 778,136\\
\code{M2} & 26,807 & 71,852 & 523,646\\
\code{F35} & 26,375 & 70,957 & 1,070,069\\
\code{F37} & 26,375 & 71,072 & 1,557,558\\
\bottomrule\end{tabular}
\caption{Uniform exclusions for $35\le N\le40$.}
\label{app:tab:certificates}
\end{table}

The total number of checked RUP additions is 7,023,863.
Two of these refutations, \code{G26} and \code{T1}, were additionally
checked using a separate Glucose 4 propagation implementation, with
163,776 checked additions in total and no calls to SAT search.
This second check concerns these two formulas, not all fifteen.

Independent implementation checks cover the pruning steps, the local
orbit partitions, and agreement of the proper-cycle data with the
constructions in Section~\ref{app:supp:cycles}. The signed weighted encoding
was tested on 6,144 assignments and the integer identities on 402,682
profiles as regression checks. Their mathematical justification is
given above and in Section~\ref{app:supp:encoding}.
A known eight-rectangle covering at $n=34$ is included as a positive
control and is accepted when the lower selection threshold is 34.
Its rectangles are checked on every matrix position.

The supplement \code{eight-rectangle-threshold-certificate.zip} contains
the fifteen CNFs and DRUP refutations, the checker and its source,
generators, audits, proper-cycle parameters, the positive control, the
published upper-bound witnesses, and a SHA-256 manifest.
Checking saved refutations does not require a SAT-solver package.
Rebuilding formulas uses \code{python-sat==1.9.dev17}.
The key mathematical lemmas are given in the main paper; the complete
parameters and encodings needed to inspect the certificates are given here.

\section{Proof checking and reproduction commands}\label{app:supp:checking}

\subsection{RUP checking and its scope}

Refutations were produced by CaDiCaL 3.0.0 through
Python-SAT 1.9.dev17~\cite{app:CaDiCaL,app:PySAT}.
A clause $L$ has the reverse-unit-propagation (RUP) property relative
to a formula $F$ if unit propagation on
$F\land\bigwedge_{\ell\in L}\neg\ell$ yields a contradiction.
Thus $F$ entails $L$. The supplied C\# checker verifies every addition
in this sense and accepts only after verifying the empty clause.
It uses no SAT-solver library. The deletion-enabled RUP format is a
restricted form of DRAT~\cite{app:DRAT}; no unchecked RAT-only step is used.

The checker may retain previously justified root-level consequences
when processing deletions. These remain entailed by the original CNF,
so retaining them preserves soundness of the refutation.
The additional Glucose check likewise retains original and verified
unit clauses while processing other deletions in batches.

Refutation checking establishes unsatisfiability of the supplied CNFs.
The local-exclusion arguments in the main paper justify their coverage of
the original mathematical problem, and Section~\ref{app:supp:cycles}
justifies their proper-cycle clauses. The trusted components include
these reductions, the encodings, the checker, and the execution
environment. No third-party audit or proof-assistant formalization is
asserted.

\subsection{Verification commands}
From the extracted archive root, the published upper witnesses can be
checked using only the Python standard library:
\begin{verbatim}
python outputs/verify_upper_covers.py outputs/upper-covers-35-55.json
\end{verbatim}
The fifteen refutations can be checked on Windows with the supplied
executable:
\begin{verbatim}
$cases = 'H43','H57','H83','H99','H113','G23','G25','G26',
         'T1','T2','T3','M0','M2','F35','F37'
foreach ($name in $cases) {
  & .\work\RupCheck2.exe "work/threshold35/$name.cnf" `
                        "work/threshold35/$name.drup"
  if ($LASTEXITCODE -ne 0) { throw "Rejected: $name" }
}
\end{verbatim}
The checker source is \code{work/RupCheck.cs}. With the pinned Python-SAT
version installed, the arithmetic audit and the formula-reconstruction
audit are
\begin{verbatim}
python work/audit_threshold_data.py
python work/audit_threshold_proofs.py
\end{verbatim}
These audits check arithmetic, controls, model reconstruction, hashes,
and recorded checker results. They do not replace running the RUP
checker on the proof files.

An unfinished search is reported as \code{UNKNOWN} and is not a proof.
The reported exclusions depend on the saved, checked refutations.

\section{Published size-nine witnesses}\label{app:supp:upper}

For completeness of the supplementary record, we include all twenty-one
published size-nine witnesses in \code{upper-covers-35-55.json}, together
with their source-file hashes and a standard-library checker
\code{verify_upper_covers.py}. Their source is the supporting data linked
from~\cite{app:GGS}; they are not new constructions of this paper.
To decode a witness, let $c_i,d_j\in\{0,1\}^9$ be its row and column
membership vectors. The nine rectangles are explicitly
\begin{equation}\label{app:eq:upper-rectangles}
 R_t=\{i:(c_i)_t=1\}\times\{j:(d_j)_t=1\},
 \qquad 0\le t\le8.
\end{equation}
The checker verifies the exact integer condition
\begin{equation}\label{app:eq:upper-check}
 \sum_{t=0}^{8}(c_i)_t(d_j)_t>0
 \quad\Longleftrightarrow\quad i\notin\{j,j+1\}\pmod n.
\end{equation}
It therefore checks both the absence of zero entries in the rectangles
and the coverage of every positive entry. All $43,295$ matrix positions
in these twenty-one witnesses pass this check.

\end{document}